\documentclass[12pt]{article}
\usepackage{amsfonts}
\usepackage{latexsym,amsmath,amssymb}
\usepackage{graphics,epstopdf}
\usepackage[pdftex]{graphicx}
\usepackage[ utf8]{inputenc}

\numberwithin{equation}{section}
\begin{document}

\bigskip

\bigskip

\begin{center}
{\Large \textbf{\ Approximation by generalized Sz\'{a}sz operators involving Sheffer polynomials}}

\bigskip

\textbf{M. Mursaleen} and \textbf{Khursheed J. Ansari}

Department of\ Mathematics, Aligarh Muslim University, Aligarh--202002, India%
\\[0pt]

mursaleenm@gmail.com; ansari.jkhursheed@gmail.com\\[0pt]

\bigskip

\bigskip

\textbf{Abstract}
\end{center}

\parindent=8mm {\footnotesize {The purpose of this article is to give a Chlodowsky type generalization of Sz\'{a}sz operators defined by means of the Sheffer type polynomials. We obtain convergence properties of our operators with the help of Korovkin's theorem and the order of convergence by using a classical approach, the second order modulus of continuity and Peetre's $K$-functional. Moreover, we study the convergence of these operators in a weighted space of functions on a positive semi-axis and estimate the approximation by using a new type of weighted modulus of continuity introduced by {\it Gadjiev and Aral} in \cite{gad3}. An algorithm is also given to plot graphical examples, and we have shown the convergence of these operators towards the function and these examples can be take as a comparison between the new operators with the previous one too. Finally, some numerical examples are also given.}\newline

{\footnotesize \emph{Keywords and phrases}: Sz\'{a}sz operators, Modulus of continuity, Rate of convergence, Weighted

\parindent8mmspace, Sheffer polynomials.}


\noindent\textbf{1. Introduction and preliminaries}\newline

\parindent=8mm In approximation theory, the positive approximation processes discovered by Korovkin play a central role and arise in
a natural way in many problems connected with functional analysis, harmonic analysis, measure theory, partial differential
equations and probability theory. The most useful examples of such operators are Sz\'{a}sz \cite{sza} operators.\\

\parindent8mm Sz\'{a}sz \cite{sza} defined the positive linear operators:
\begin{equation*}
    S_n(f;x):=e^{-nx}\sum\limits_{k=0}^{\infty}\frac{(nx)^k}{k!}f\bigg(\frac kn\bigg)\eqno(1.1)
\end{equation*}
where $x\geq0$ and $f\in C[0,\infty)$ whenever the above sum converges. Motivated by this work, many authors have investigated several interesting properties of the operators (1.1).

\parindent8mm Later, Jakimovski and Leviatan \cite{ibr} obtained a generalization of Sz\'{a}sz operators by means of Appell polynomials. Let $g(z)=\sum_{k=0}^{\infty}a_kz^k~(a_0\ne0)$ be an analytic function in the disk $|z|<R,~(R>1)$ and suppose that $g(1)\ne0$. The Appell polynomials $p_k(x)$ have generating functions of the form
\begin{equation*}
    g(u)e^{ux}=\sum\limits_{k=0}^{\infty}p_k(x)u^k.\eqno(1.2)
\end{equation*}
Under the assumption that $p_k(x)\geq0$ for $x\in[0,\infty)$, Jakimovski and Leviatan introduced the positive linear operators $P_n(f;x)$ via

\begin{equation*}
    P_n(f;x):=\frac{e^{-nx}}{g(1)}\sum\limits_{k=0}^{\infty}p_k(nx)f\bigg(\frac kn\bigg)\eqno(1.3)
\end{equation*}
and gave the approximation properties of the operators.\\

{\it Case} 1. For $g(1)=1$, with the help of (1.2) we easily find $p_k(x)=\frac{x^k}{k!}$ and from (1.3) we meet again the Sz\'{a}sz operators given by (1.1).

\vspace{.25cm}

\parindent8mm Then, Ismail \cite{ism} presented another generalization of Sz\'{a}sz operators (1.1) and Jakimovski and Leviatan operators (1.3) by using Sheffer polynomials. Let $A(z)=\sum_{k=0}^{\infty}a_kz^k~(a_0\ne0)$ and $H(z)=\sum_{k=1}^{\infty}h_kz^k~(h_1\ne0)$ be analytic functions in the disk $|z|<R~(R>1)$ where $a_k$ and $h_k$ are real. The Sheffer polynomials $p_k(x)$ have generating functions of the type

\begin{equation*}
    A(t)e^{xH(t)}=\sum\limits_{k=0}^{\infty}p_k(x)t^k,~~~|t|<R.\eqno(1.4)
\end{equation*}
Using the following assumptions:
\begin{itemize}
  \item [(i)] for $x\in[0,\infty)$,  $p_k(x)\geq0$,
  \item [(ii)] $A(1)\ne0$ and $H^\prime(1)=1$,\hspace{10cm}(1.5)
\end{itemize}

\noindent Ismail investigated the approximation properties of the positive linear operators given by
\begin{equation*}
    T_n(f;x):=\frac{e^{-nxH(1)}}{A(1)}\sum\limits_{k=0}^{\infty}p_k(nx)f\bigg(\frac kn\bigg), \mbox{  for }n\in\mathbb{N}.\eqno(1.6)
\end{equation*}
{\it Case} 1. For $H(t)=t$, it can be easily seen that the generating functions (1.4) return to (1.2) and, from this fact, the operators (1.6) reduce to the operators (1.3).\\
\\
{\it Case} 2. For $H(t)=t$ and $A(t)=1$, one get the Sz\'{a}sz operators from the operators (1.6).\newline

\parindent8mm In \cite{ibr}, B\"{u}y\"{u}kyaz{\i}c{\i} et al. introduced the Chlodowsky \cite{chl} variant of operators (1.3). Guided by their work we give the Chlodowsky type generalization of operators (1.6) as follows:

\begin{equation*}
    T_n^*(f;x):=\frac{e^{-\frac{n}{b_n}xH(1)}}{A(1)}\sum\limits_{k=0}^{\infty}p_k\bigg(\frac{n}{b_n}x\bigg)f\bigg(\frac knb_n\bigg)\eqno(1.7)
\end{equation*}

\noindent with $b_n$ a positive increasing sequence with the properties
\begin{equation*}
    \lim\limits_{n\to\infty}b_n=\infty,~~~~\lim\limits_{n\to\infty}\frac{b_n}{n}=0\eqno(1.8)
\end{equation*}
and $p_k$ are Sheffer polynomials defined by (1.4). For other generalization of operators (1.6) one can refer to \cite{sez}. \\

\parindent 8mm The rest of the paper is organized as follows. In Section 2 we obtain some local approximation results by the generalized Sz\'{a}sz operators given by (1.7). In particular, the convergence of operators is examined with the help of Korovkin's theorem. The order of approximation is established by means of a classical approach, the second-order modulus of continuity and Peetre's $K$-functional. An algorithm and some graphical examples are also given to in claim of convergence of operators towards the function. Section 3 is devoted to study the convergence of these operators in a weighted space of functions on a positive semi-axis and estimate the approximation by using a new type of weighted modulus of continuity introduced by {\it Gadjiev and Aral} in \cite{gad3}. Finally, some numerical examples are also given in section 4.\\

\parindent8mm Note that throughout the paper we will assume that the operators $T_n^*$ are positive and we use the following test functions
\begin{equation*}
    e_i(x)=x^i,~~~~i\in\{0,1,2\}.
\end{equation*}

\vspace{.5cm}

\noindent\textbf{2. Local approximation properties of $T_n^*(f;x)$}%
\newline

\parindent=8mm We denote by $C_E[o,\infty)$ the set of all continuous functions $f$ on $[0,\infty)$ with the property that $|f(x)|\leq\beta e^{\alpha x}$ for all $x\geq0$ and some positive finite $\alpha$ and $\beta$. For a fixed $r\in\mathbb{N}$ we denote by $C_E^r[0,\infty)=\{f\in C_E[0,\infty):f^\prime,f^{\prime\prime},\cdots,f^{(r)}\in C_E[0,\infty)\}$. Using equality (1.1) and the fundamental properties of the $T_n^*$ operators, one can easily get the following lemmas:

\vspace{.5cm}

\noindent\textbf{Lemma 2.1.} For all $x\in[0,\infty)$, we have
\begin{eqnarray*}
T_n^*(e_0;x)&=& 1;\hspace{12cm}(2.1)\\
T_n^*(e_1;x)&=& x+\frac{b_n}{n}\frac{A^\prime(1)}{A(1)};\hspace{10.2cm}(2.2)\\
T_n^*(e_2;x)&=& x^2+\frac{b_n}{n}\frac{A(1)
+2A^\prime(1)+A(1)H^{\prime\prime}(1)}{A(1)}x
+\frac{b_n^2}{n^2}\frac{A^\prime(1)+A^{\prime\prime}(1)}{A(1)}.\hspace{3.5cm}(2.3)\\
\end{eqnarray*}

\parindent8mm It follows from Lemma 2.1 that,
\begin{equation*}
    T_n^*\big((e_1-x);x\big)= 0,\eqno(2.4)
\end{equation*}

\begin{equation*}
T_n^*\big((e_1-x)^2;x\big)=
\frac{b_n}{n}\big(1+H^{\prime\prime}(1)\big)x
+\frac{b_n^2}{n^2}\frac{A^\prime(1)+A^{\prime\prime}(1)}{A(1)}.\eqno(2.5)
\end{equation*}

\vspace{.5cm}

\noindent\textbf{Theorem 2.2.} For $f\in C_E[0,\infty)$, the operators $T_n^*$ converge uniformly to $f$ on $[0,a]$ as $n\to\infty$.\\
\\
\textbf{Proof.} According to (2.1)-(2.3), we have
\begin{equation*}
    \lim\limits_{n\to\infty}T_n^*(e_i;x)=e_i(x),i\in\{0,1,2\}.
\end{equation*}
If we apply the Korovkin theorem \cite{alt}, we obtain the desired result.

\vspace{.5cm}

\noindent\textbf{Algorithm}\newline

\noindent Graphically, to show the approximation of a given function $f(x)$ by positive linear operators $\mathbf{T}_n(f;x)$ and $\mathbf{T}_n^*(f;x)$ given by (1.7), the algorithm is summarized as below.\newline

\noindent\textbf{Step 1:} Choose the functions $A(t)$ and $H(t)$ such that $A(t)\ne0$ and $H^\prime(1)=1$.\newline
\noindent\textbf{Step 2:} Find out the Sheffer polynomials $\mathbf{p}_k(x)$ with the help of relation (1.4).\newline
\noindent\textbf{Step 3:} Check $\mathbf{p}_k(x)\ge0$ for $x\ge0$.\newline
\noindent\textbf{Step 4:} Choose the sequence $b_n$ under the condition given in (1.8).\newline
\noindent\textbf{Step 5:} Plot the graph of function $f(x)$ and the operators $\mathbf{T}_n(f;x)$ and $\mathbf{T}_n^*(f;x)$ for the different values of $n$.\newline

\vspace{.5cm}

\noindent\textbf{Example 2.3.} For (i) $A(t)=e^t \mbox{ and } H(t)=t$, (ii) $A(t)=t \mbox{ and } H(t)=t$, the convergence of the two operators $T_n(f;x)$ and $T_n^*(f;x)$ to $f(x)$ are illustrated in Figs. 1, 2, 3 and 4 respectively, where $f(x)=-4xe^{-3x},~n=10,~50,~100,~200,~300$, and $b_n=\sqrt{n}$.

\begin{figure*}[htb!]
\begin{center}
\includegraphics[height=8cm, width=15cm]{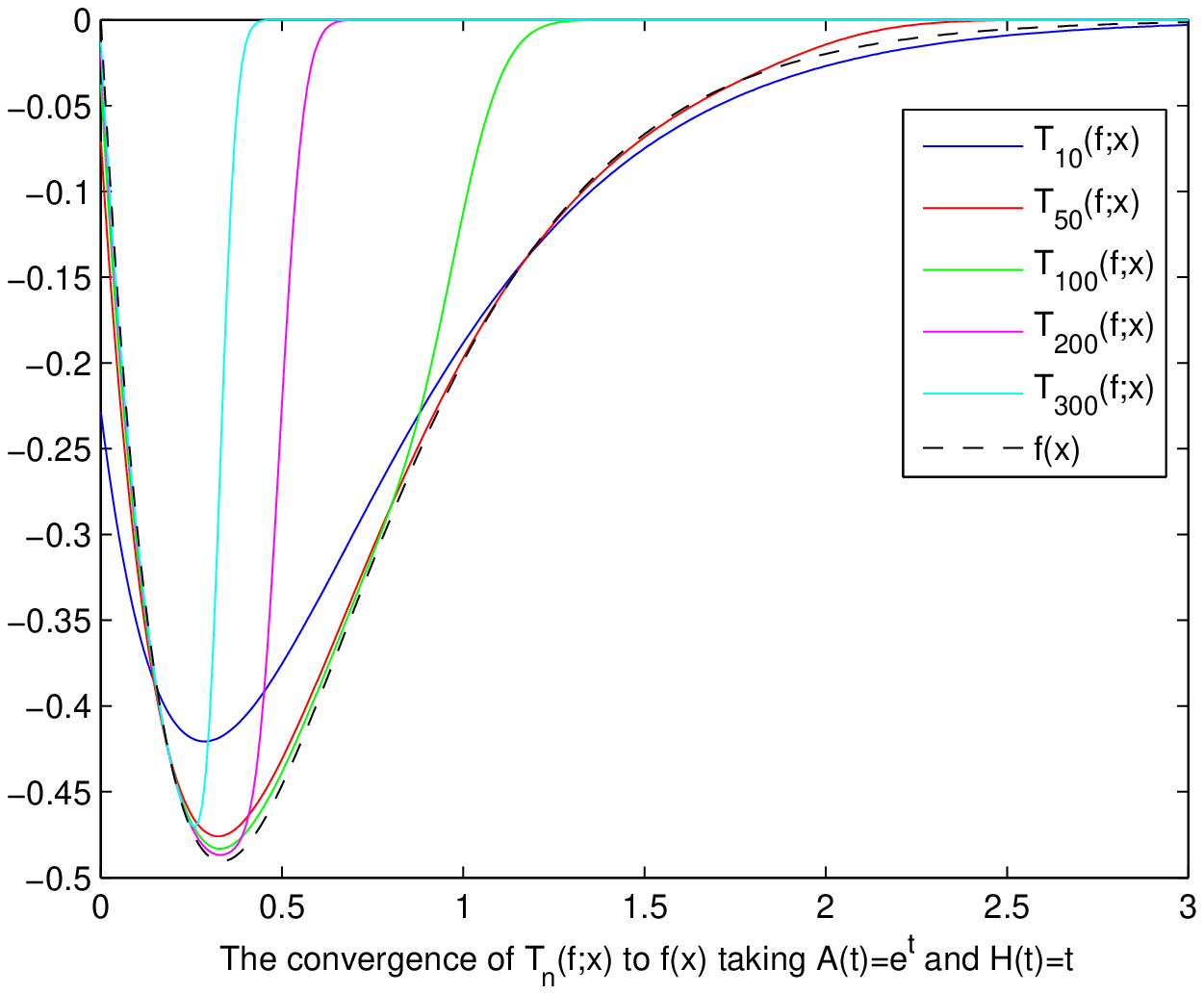}
\end{center}
\caption{}
\end{figure*}

\begin{figure*}[htb!]
\begin{center}
\includegraphics[height=8cm, width=15cm]{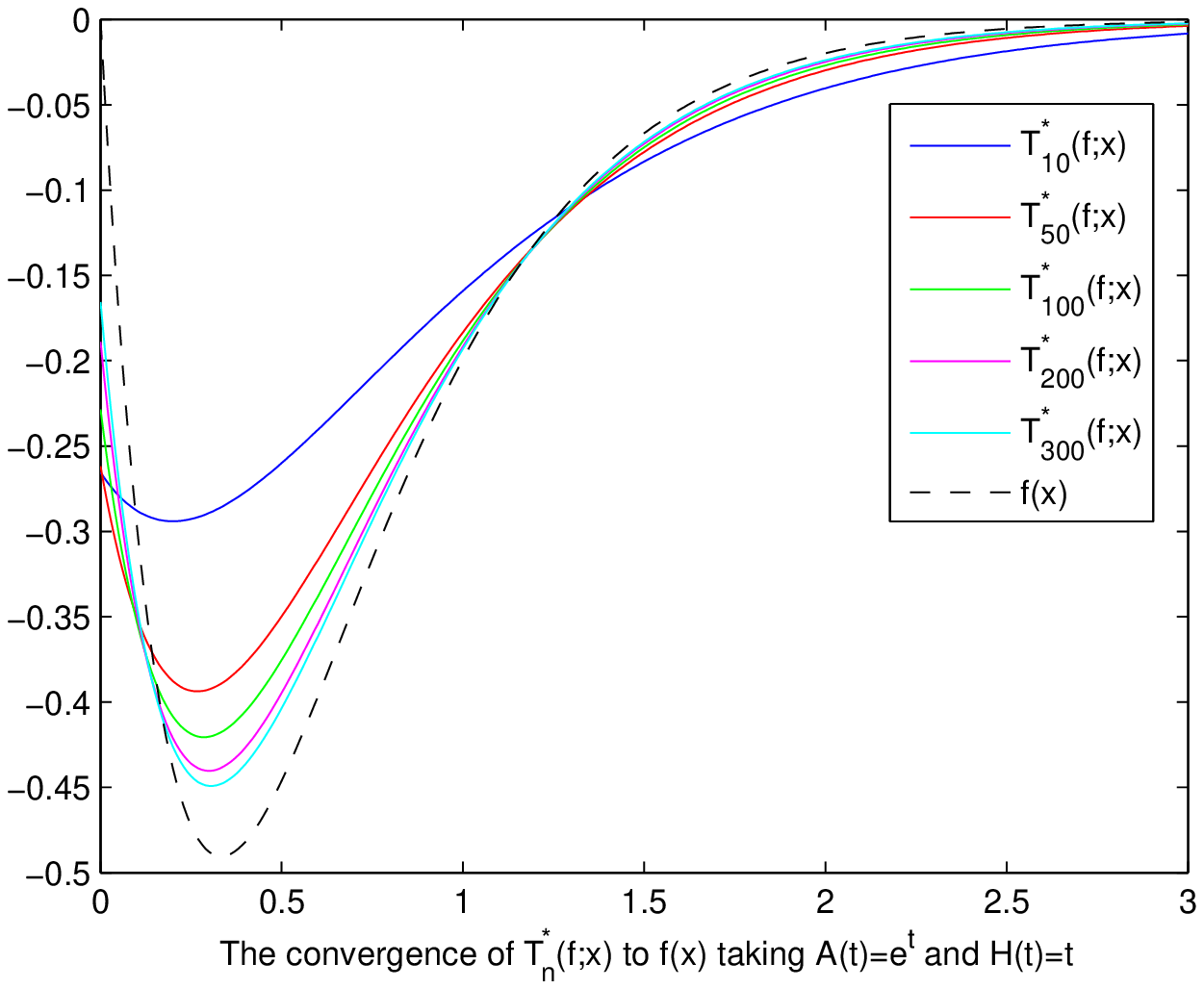}
\end{center}
\caption{}
\end{figure*}

\begin{figure*}[htb!]
\begin{center}
\includegraphics[height=8cm, width=15cm]{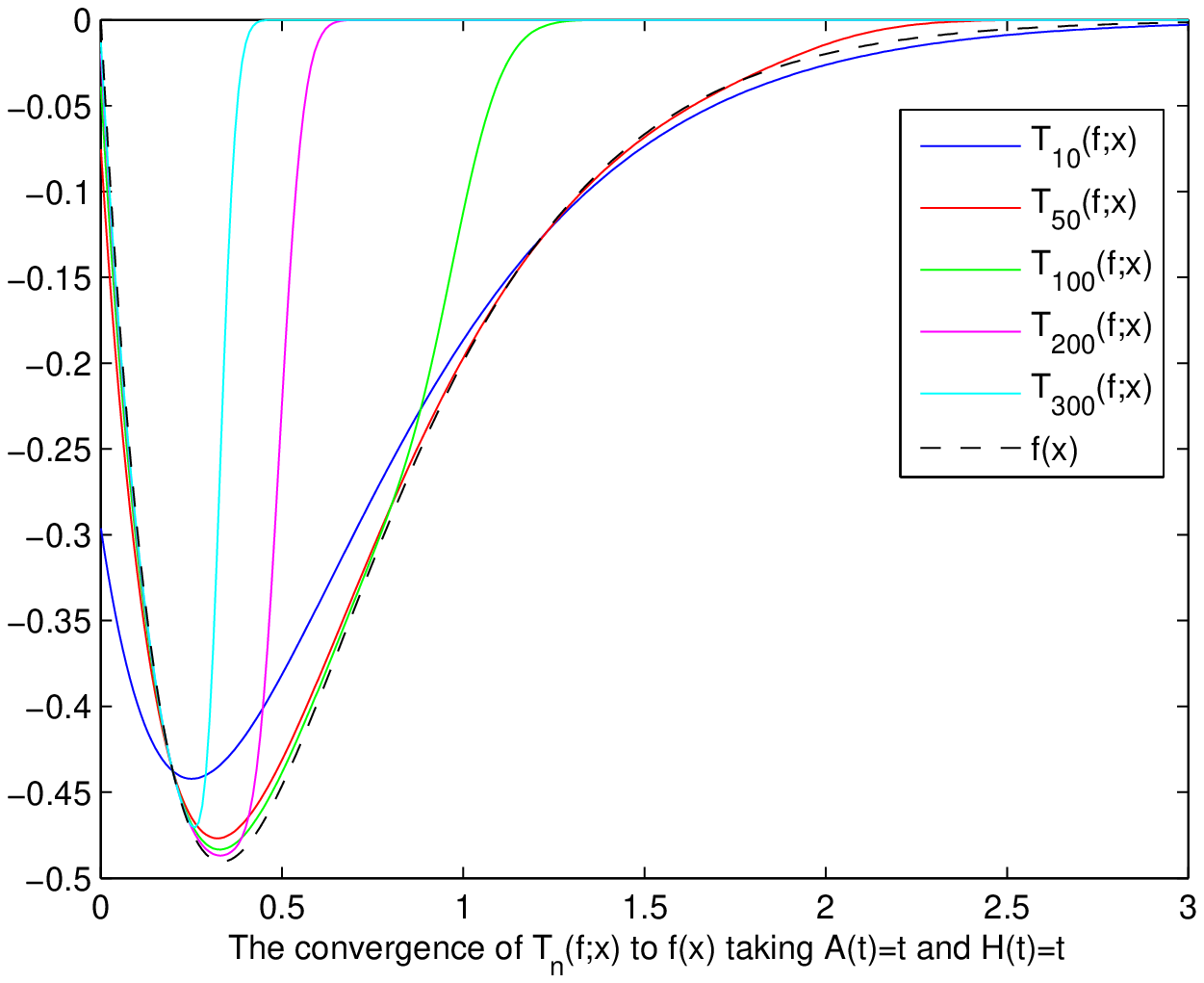}
\end{center}
\caption{}
\end{figure*}

\begin{figure*}[htb!]
\begin{center}
\includegraphics[height=8cm, width=15cm]{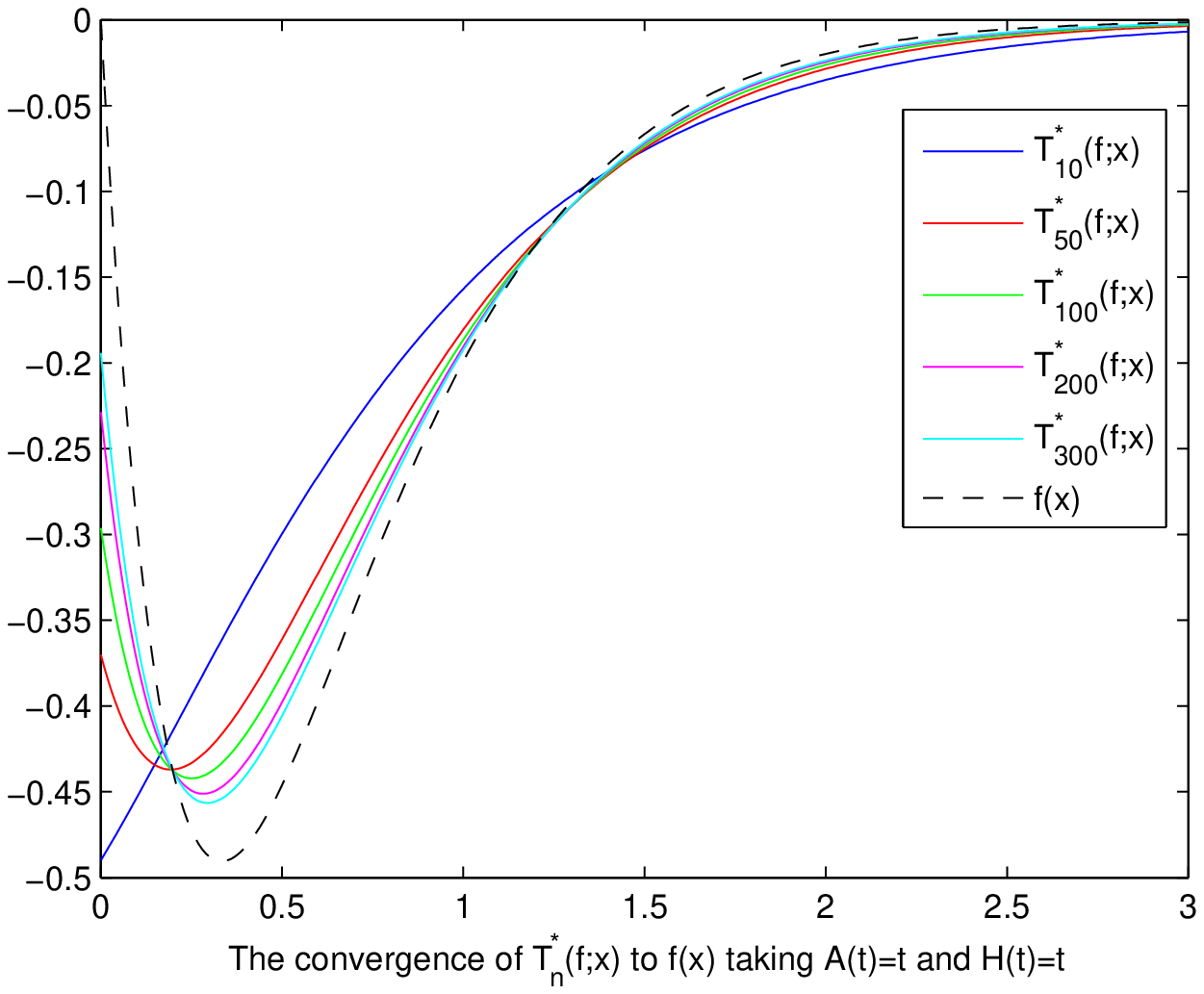}
\end{center}
\caption{}
\end{figure*}

\vspace{.5cm}
\noindent\textbf{Remark.}
 From figure 1 and figure 3, we can see that when the value of $n$ is increasing, the graph of operators $T_n(f;x)$ are going far away form the graph of the function $f(x)$ but with our proposed operators $T_n^*(f;x)$, the convergence towards the function can be seen very clearly from figure 2 and figure 4. In a nut shell, we can claim that to approximate a function our operators $T_n^*(f;x)$ are better in comparison of $T_n(f;x)$.

\vspace{.5cm}
\parindent8mm Now, we concerned with the estimate of the order of approximation of a function $f$ by means of the positive operators $L_n^*$, using the first and second order modulus of continuity \cite{dit}.\\

\vspace{.25cm}
\parindent8mm If $\delta>0$, the modulus of continuity $\omega(f,\delta)$ of $f\in C[a,b]$ is defined by
\begin{equation*}
    \omega(f,\delta)=\sup\limits_{x,y\in[a,b],~|x-y|\leq\delta}|f(x)-f(y)|.
\end{equation*}

\vspace{.25cm}
\parindent8mm The second order modulus of continuity of $f\in C_B[0,\infty)$ is defined by
\begin{equation*}
    \omega_2(f,\delta)=\sup\limits_{0<t\leq\delta}\|f(.+2t)-2f(.+t)+f(.)\|_{C_B}
\end{equation*}
where $C_B[0,\infty)$ is the class of real valued functions defined on $[0,\infty)$ which are bounded and uniformly continuous with the norm $\|f\|_{C_B}=\sup_{x\in[0,\infty)}|f(x)|$.\\

\vspace{.25cm}
\parindent8mm The Peetre's $K$-functional \cite{dit} of the function $f\in C_B[0,\infty)$ is defined by
\begin{equation*}
    K(f,\delta):=\inf\limits_{g\in C_B^2[0,\infty)}\big\{\|f-g\|_{C_B}+\delta\|g\|_{C_B^2}\big\}
\end{equation*}
where
\begin{equation*}
    C_B^2[0,\infty):=\big\{g\in C_B[0,\infty):~g^\prime,~g^{\prime\prime}\in C_B[0,\infty)\big\}
\end{equation*}
and the norm $\|g\|_{C_B^2}:=\|g\|_{C_B}+\|g^\prime\|_{C_B}+\|g^{\prime\prime}\|_{C_B}$. It is clear that the following inequality:
\begin{equation*}
    K(f,\delta)\leq M\big\{\omega_2(f,\sqrt{\delta})+\min(1,\delta)\|f\|_{C_B}\big\}
\end{equation*}
is valid, for all $\delta>0$. The constant $M$ is independent of $f$ and $\delta$.\\

\vspace{.25cm}

\noindent\textbf{Lemma 2.4} {\it(\cite{gav})}{\bf.} Let $g\in C^2[0,\infty)$ and $(P_n)_{n\geq0}$ be a sequence of positive linear operators with the property $P_n(1;x)=1$. Then
\begin{equation*}
    |P_n(g;x)-g(x)|\leq\|g^{\prime}\|\sqrt{P_n\big((s-x)^2;x\big)}+\frac12 \|g^{\prime\prime}\|P_n\big((s-x)^2;x\big).
\end{equation*}

\vspace{0.25cm}

\noindent\textbf{Lemma 2.5} {\it(\cite{zhu})}{\bf.} Let $f\in C[a,b]$ and $h\in\big(0,\frac{b-a}2\big)$. Let $f_h$ be the second-order Steklov function attached to the function $f$. Then the following inequalities are satisfied:
\begin{itemize}
  \item [(i)] $\|f_h-f\|\leq\frac34\omega_2(f,h)$,
  \item [(ii)] $\|f_h^{\prime\prime}\|\leq\frac{3}{2h^2}\omega_2(f,h)$.
\end{itemize}

\vspace{.25cm}

\noindent\textbf{Theorem 2.6.} If $f\in C_E[0,\infty)$, then for any $x\in[0,a]$ we have
\begin{equation*}
    |T_n^*(f;x)-f(x)|\leq\Bigg\{1+\sqrt{\big(1+H^{\prime\prime}(1)\big)a+\frac{b_n}{n}\frac{A^\prime(1)+A^{\prime\prime}}{A(1)}}\Bigg\}
    \omega\Bigg(f,\sqrt{\frac{b_n}{n}}\Bigg).
\end{equation*}

\vspace{.25cm}
\noindent\textbf{Proof.} We will use the relation (2.5) and the well-known properties of the modulus of continuity. We have
\begin{eqnarray*}
  |T_n^*(f;x)-f(x)|&\leq& \frac{e^{-\frac{n}{b_n}xH(1)}}{A(1)}\sum\limits_{k=0}^{\infty}p_k\Big(\frac{n}{b_n}x\Big)\Big|f\Big(\frac knb_n\Big)-f(x)\Big|\\
  &\leq&\bigg\{1+\frac1\delta\frac{e^{-\frac{n}{b_n}xH(1)}}{A(1)}\sum\limits_{k=0}^{\infty}p_k\Big(\frac{n}{b_n}x\Big)\Big|\frac knb_n-x\Big|\bigg\}\omega(f,\delta).
\end{eqnarray*}
Recalling the Cauchy-Schwartz inequality we obtain the formula below:
\begin{eqnarray*}
  |T_n^*(f;x)-f(x)|  &\leq&\bigg\{1+\frac1\delta\bigg(\frac{e^{-\frac{n}{b_n}xH(1)}}{A(1)}\sum\limits_{k=0}^{\infty}p_k\Big(\frac{n}{b_n}x\Big)\Big(\frac knb_n-x\Big)^2\bigg)^{\frac12}\bigg\}\omega(f,\delta)\\
  &=&\Big\{1+\frac1\delta\sqrt{T_n^*\big((e_1-x)^2;x\big)}\Big\}\omega(f,\delta).\hspace{6cm}(2.6)
\end{eqnarray*}
By means of (2.5), for $0\leq x\leq a$, one gets
\begin{equation*}
T_n^*\big((t-x)^2;x\big)\leq
\frac{b_n}{n}\big(1+H^{\prime\prime}(1)\big)a
+\frac{b_n^2}{n^2}\frac{A^\prime(1)+A^{\prime\prime}(1)}{A(1)}.\eqno(2.7)
\end{equation*}
Using (2.7) and taking $\delta=\sqrt{\frac{b_n}{n}}$ in (2.9), we obtain the desired result.\\

\vspace{.25cm}

\noindent\textbf{Theorem 2.7.} For $f\in C[0,a]$, the following inequality:
\begin{equation*}
    |T_n^*(f;x)-f(x)|\leq\frac 2a\|f\|h^2+\frac34(a+2+h^2)\omega_2(f,h)
\end{equation*}
is satisfied where
\begin{equation*}
    h:=h_n(x)=\sqrt[4]{T_n^*\big((e_1-x)^2;x\big)}
\end{equation*}
and the second order modulus of continuity is given by $\omega_2(f,\delta)$ with the norm $\|f\|=\max\limits_{x\in[a,b]}|f(x)|$.\\
\\
\textbf{Proof.} Let $f_h$ be the second-order Steklov function attached to the function $f$. By virtue of the identity (2.1), we have
\begin{eqnarray*}
  |T_n^*(f;x)-f(x)|&\leq&|T_n^*(f-f_h;x)|+|T_n^*(f_h;x)-f_h(x)|+|f_h(x)-f(x)| \\
  &\leq&2\|f_h-f\|+|T_n^*(f_h;x)-f_h(x)|.\hspace{5.4cm}(2.8)
\end{eqnarray*}
Taking into account the fact that $f_h\in C^2[0,a]$, it follows from Lemma 2.4 that
\begin{equation*}
    |T_n^*(f_h;x)-f_h(x)|\leq\|f_h^{\prime}\|\sqrt{T_n^*\big((e_1-x)^2;x\big)}+\frac12 \|f_h^{\prime\prime}\|T_n^*\big((e_1-x)^2;x\big).\hspace{3cm}(2.9)
\end{equation*}
Combining the Landau inequality and Lemma 2.5, we can write
\begin{eqnarray*}
      \|f_h^\prime\|&\leq&\frac2a\|f_h\|+\frac a2\|f_h^{\prime\prime}\|\\
      &\leq&\frac2a\|f\|+\frac{3a}{4}\frac{1}{h^2}\omega_2(f,h).
\end{eqnarray*}
From the last inequality, (2.9) becomes, on taking $h=\sqrt[4]{T_n^*\big((e_1-x)^2;x\big)}$,
\begin{equation*}
    |T_n^*(f_h;x)-f_h(x)|\leq\frac2a\|f\|h^2+\frac{3a}{4}\omega_2(f,h)+\frac{3}{4}h^2\omega_2(f,h).\hspace{.8cm}\eqno(2.10)
\end{equation*}
Substituting (2.10) in (2.9), Lemma 2.5 hence gives the proof of the theorem.\\

\vspace{.25cm}

\noindent\textbf{Theorem 2.8.} Let $f\in C_B^2[0,\infty)$. Then
\begin{equation*}
    |T_n^*(f;x)-f(x)|\leq \gamma_n(x)\|f\|_{C_B^2}
\end{equation*}
where
\begin{eqnarray*}
 \gamma(x):=\gamma_n(x)=\frac12~T_n^*\big((t-x)^2;x\big).
\end{eqnarray*}
\\
\textbf{Proof.} Using the Taylor expansion of $f$, the linearity property of the operators $T_n^*$ and (2.1), it follows that
\begin{equation*}
    T_n^*(f;x)-f(x)=f^\prime(x)T_n^*(e_1-x;x)+\frac12f^{\prime\prime}(\eta)T_n^*\big((e_1-x)^2;x\big),~~~~\eta\in(x,t).\eqno(2.11)
\end{equation*}
Taking into account the fact that
\begin{equation*}
    T_n^*\big((e_1-x);x\big)\geq0
\end{equation*}
for $x\leq t$, by combining Lemmas 2.1 and (2.5) in (2.11) we are led to
\begin{eqnarray*}
  T_n^*(f;x)-f(x)&\leq& \frac12~\big\{T_n^*\big((t-x)^2;x\big)\big\}\|f^{\prime\prime}\|_{C_B}\\
   &\leq&  \frac12~\big\{T_n^*\big((t-x)^2;x\big)\big\}\|f\|_{C_B^2}
\end{eqnarray*}

\noindent which completes the proof.\\

\vspace{.25cm}
\noindent\textbf{Theorem 2.9.} Let $f\in C_B[0,\infty)$. Then
\begin{equation*}
|T_n^*(f;x)-f(x)|\leq2M\big\{\omega_2(f,\sqrt{\delta})+\min(1,\delta)\|f\|_{C_B}\big\}
\end{equation*}
where
$\delta:=\delta_n(x)=\frac14\gamma_n(x)$ and $M>0$ is a constant independent of the function $f$ and $\delta$. Note that $\gamma_n(x)$ is defined as in Theorem 2.8.\\
\\
\textbf{Proof.} Let $g\in C_B^2[0,\infty)$. Theorem 2.8 allows us to write
\begin{eqnarray*}
  |T_n^*(f;x)-f(x)|&\leq&|T_n^*(f-g;x)|+|T_n^*(g;x)-g(x)|+|g(x)-f(x)|  \\
   &\leq&2\|f-g\|_{C_B}+\frac12~\big\{T_n^*\big((t-x)^2;x\big)\big\}\|g\|_{C_B^2}\\
   &=&2\big\{\|f-g\|_{C_B}+\delta\|g\|_{C_B^2}\big\}.\hspace{7cm}(2.12)
\end{eqnarray*}
The left-hand side of inequality (2.12) does not depend on the function $g\in C_B^2[0,\infty)$, so
\begin{equation*}
     |T_n^*(f;x)-f(x)|\leq2K(f,\delta).\eqno(2.13)
\end{equation*}
By using the relation between Peetre's $K$-functional and second modulus of smoothness, (2.13) becomes
\begin{equation*}
    |T_n^*(f;x)-f(x)|\leq2M\big\{\omega_2(f,\sqrt{\delta})+\min(1,\delta)\|f\|_{C_B}\big\}.
\end{equation*}

\vspace{.5cm}

\noindent\textbf{3. Approximation properties in weighted spaces}%
\\

\parindent8mm Now we give approximation properties of the operators $T_n^*$ of the weighted spaces of continuous functions with exponential growth on $\mathbb{R}_0^+=[0,\infty)$ with the help of the weighted Korovkin type theorem proved by Gadjiev in \cite{gad1,gad2}. For this purpose, we consider the following weighted spaces of functions which are defined on the $\mathbb{R}_0^+=[0,\infty)$.

\vspace{.25cm}
\parindent8mm Let $\rho(x)$ be the weighted function and $M_f$ a positive constant, then we define
\begin{eqnarray*}
   B_{\rho}(\mathbb{R}_0^+)&=& \{f\in E(\mathbb{R}_0^+): |f(x)\leq M_f\rho(x)|\},\\
   C_{\rho}(\mathbb{R}_0^+)&=& \{f\in B_{\rho}(\mathbb{R}_0^+): f \mbox{~is continuous}\},\\
   C_{\rho}^k(\mathbb{R}_0^+)&=& \bigg\{f\in C_{\rho}(\mathbb{R}_0^+): \lim\limits_{n\to\infty}\frac{f(x)}{\rho(x)}=K_f<\infty\bigg\}.
\end{eqnarray*}
It is obvious that $C_{\rho}^k(\mathbb{R}_0^+)\subset C_{\rho}(\mathbb{R}_0^+)\subset B_{\rho}(\mathbb{R}_0^+)$. The space $B_{\rho}(\mathbb{R}_0^+)$ is a normed linear space with the following norm:
\begin{equation*}
    \|f\|_{\rho}=\sup\limits_{x\in \mathbb{R}_0^+}\frac{|f(x)|}{\rho(x)}.
\end{equation*}
The following results on the sequence of positive linear operators in these spaces are given \cite{gad1,gad2}.\\

\vspace{.25cm}
\noindent\textbf{Lemma 3.1} {\it(\cite{gad1,gad2})}\textbf{.} The sequence of positive linear operators $(L_n)_{n\geq1}$ which act from $C_{\rho}(\mathbb{R}_0^+)$ to $B_{\rho}(\mathbb{R}_0^+)$ if and only if there exists a positive constant $k$ such that
\begin{equation*}
    L_n(\rho;x) \leq k\rho(x),~~~\mbox{i.e.}
\end{equation*}
\begin{equation*}
    \|L_n(\rho;x)\|_{\rho}\leq k.
\end{equation*}

\vspace{.25cm}
\noindent\textbf{Theorem 3.2} {\it(\cite{gad1,gad2})}\textbf{.} Let $(L_n)_{n\geq1}$ be the sequence of positive linear operators which act from $C_{\rho}(\mathbb{R}_0^+)$ to $B_{\rho}(\mathbb{R}_0^+)$ satisfying the conditions
\begin{equation*}
    \lim\limits_{n\to\infty}\|L_n(t^i;x)-x^i\|_{\rho}=0,~~~i\in\{0,1,2\},
\end{equation*}
then for any function $f\in C_{\rho}^k(\mathbb{R}_0^+)$
\begin{equation*}
    \lim\limits_{n\to\infty}\|L_nf-f\|_{\rho}=0.
\end{equation*}

\vspace{.25cm}
\noindent\textbf{Lemma 3.3.} Let $\rho(x)=1+x^2$ be a weight function. If $f\in C_{\rho}(\mathbb{R}_0^+)$, then
\begin{equation*}
    \|T_n^*(\rho;x)\|_{\rho}\leq1+M.
\end{equation*}

\noindent\textbf{Proof.} Using (2.1) and (2.3), one has
\begin{eqnarray*}
  T_n^*(\rho;x) &=& 1+x^2+\frac{b_n}{n}\frac{A(1)
+2A^\prime(1)+A(1)H^{\prime\prime}(1)}{A(1)}x
+\frac{b_n^2}{n^2}\frac{A^\prime(1)+A^{\prime\prime}(1)}{A(1)}\\
   \|T_n^*(\rho;x)\|_{\rho}&=& \sup\limits_{x\geq0}\bigg\{\frac1{1+x^2}\bigg(1+x^2+\frac{b_n}{n}\frac{A(1)
+2A^\prime(1)+A(1)H^{\prime\prime}(1)}{A(1)}x
+\frac{b_n^2}{n^2}\frac{A^\prime(1)+A^{\prime\prime}(1)}{A(1)}\bigg)\bigg\}\\
&\leq&1+\frac{b_n}{n}\frac{A(1)
+2A^\prime(1)+A(1)H^{\prime\prime}(1)}{A(1)}
+\frac{b_n^2}{n^2}\frac{A^\prime(1)+A^{\prime\prime}(1)}{A(1)}.
\end{eqnarray*}
Since $\lim_{n\to\infty}\frac{b_n}{n}=0$, there exists a positive $M$ such that
\begin{equation*}
    \|T_n^*(\rho;x)\|_{\rho}\leq1+M
\end{equation*}
so the proof is completed.\\

\parindent8mm By using Lemma 3.3, we can easily see that the operators $T_n^*$ defined by (1.7) act from $C_{\rho}(\mathbb{R}_0^+)$ to $B_{\rho}(\mathbb{R}_0^+)$.\\

\vspace{.25cm}
\noindent\textbf{Theorem 3.4.} Let $T_n^*$ be the sequence of positive linear operators defined by (1.7) and $\rho(x)=1+x^2$, then for each $f\in C_{\rho}^k(\mathbb{R}_0^+)$
\begin{equation*}
    \lim\limits_{n\to\infty}\|T_n^*(f;x)-f(x)\|_{\rho}=0.
\end{equation*}

\noindent\textbf{Proof.} It is enough to prove that the conditions of the weighted Korovkin type theorem given by Theorem 3.2 are satisfied. From (2.1), it is immediately seen that
\begin{equation*}
    \lim\limits_{n\to\infty}\|T_n^*(e_0;x)-e_0(x)\|_{\rho}=0.\eqno(3.1)
\end{equation*}
Using (2.2) we have
\begin{equation*}
\|T_n^*(e_1;x)-e_1(x)\|_{\rho}=\frac{b_n}{n}\frac{A^\prime(1)}{(1)}\eqno(3.2)
\end{equation*}
this implies that
\begin{equation*}
    \lim\limits_{n\to\infty}\|T_n^*(e_1;x)-e_1(x)\|_{\rho}=0.\eqno(3.3)
\end{equation*}
By means of (2.3) we get
\begin{eqnarray*}
  \|T_n^*(e_2;x)-e_2(x)\|_{\rho}&=&\sup\limits_{x\in R_0}\bigg|\frac{b_n}{n}\frac{A(1)
+2A^\prime(1)+A(1)H^{\prime\prime}(1)}{A(1)}\frac{x}{1+x^2}
+\frac{b_n^2}{n^2}\frac{A^\prime(1)+A^{\prime\prime}(1)}{A(1)}\frac{1}{1+x^2}\bigg|\\
&\leq&\frac{b_n}{n}\frac{A(1)
+2A^\prime(1)+A(1)H^{\prime\prime}(1)}{A(1)}+\frac{b_n^2}{n^2}\frac{A^\prime(1)+A^{\prime\prime}(1)}{A(1)}.\hspace{2.5cm}(3.4)
\end{eqnarray*}
Using the conditions (1.8), it follows that
\begin{equation*}
    \lim\limits_{n\to\infty}\|T_n^*(e_2;x)-e_2(x)\|_{\rho}=0.\eqno(3.5)
\end{equation*}
From (3.1), (3.2) and (3.5), for $i\in\{0,1,2\}$ we have
\begin{equation*}
    \lim\limits_{n\to\infty}\|T_n^*(e_i;x)-e_i(x)\|_{\rho}=0.
\end{equation*}
If we apply Theorem 3.2, we obtain the desired result.

\vspace{.25cm}

\parindent8mm Now, for any weighted function $\rho(x)$, we want to find the approximation and rate of approximation of the functions $f\in C^k_\rho(\mathbb{R}^+_0)$ by using the operators $T_n^*$ on $\mathbb{R}^+_0=[0,\infty)$. For this we define new positive linear operators which are a generalization of the $T_n^*$ operators. It is well-known that the usual first modulus of continuity does not tend to zero as $\delta\to0$ on $\mathbb{R}^+_0$, so we use the following new type of weighted modulus of continuity introduced by Gadjiev and Aral in \cite{gad3}:
\begin{equation*}
 \Omega_\rho(f,\delta)=\Omega(f,\delta)_{\mathbb{R}_0^+}
 =\sup\limits_{x,t\in\mathbb{R}_0^+,~|\rho(t)-\rho(x)|\leq\delta}\frac{|f(t)-f(x)|}{[|\rho(t)-\rho(x)|+1]\rho(x)}\eqno(3.6)
\end{equation*}
where $\rho$ is satisfying the following assumptions:
\begin{itemize}
  \item [(i)] $\rho$ is a continuously differentiable function on $\mathbb{R}_0^+$ and $\rho(0)=1$,
  \item [(ii)] $\inf_{x\geq0}\rho^\prime(x)\geq1$.
\end{itemize}

\parindent8mm The weighted modulus of continuity $\Omega_\rho(f,\delta)$ given by (3.6) has some properties given in the following lemma (see \cite{gad3}).

\vspace{.5cm}
\noindent\textbf{Lemma 3.5} {\it \cite{gad3}} {\bf.} For any $f\in C^k_\rho(\mathbb{R}_0^+)$ then
\begin{equation*}
    \lim\limits_{\delta\to0}\Omega_\rho(f,\delta)=0,
\end{equation*}
and for each $x,t\in\mathbb{R}_0^+$ the inequality
\begin{equation*}
    |f(t)-f(x)|\leq2\rho(x)(1+\delta^2)\bigg(1+\frac{\big(\rho(t)-\rho(x)\big)^2}{\delta^2}\bigg)\Omega_\rho(f,\delta)
\end{equation*}
holds, where $\delta$ is any fixed positive number.

\vspace{.25cm}

\parindent8mm The estimates of the approximation of functions by positive linear operators by means of the new type of modulus of continuity are given in the following theorem \cite{gad3}:

\vspace{.5cm}
\noindent\textbf{Theorem 3.6} {\it \cite{gad3}} {\bf.} Let $\rho(x)\leq\psi_k(x),~k=1,2,3$ and the sequences of the positive linear operators $(L_n)_{n\geq1}$ satisfying the conditions
\begin{equation*}
    \|L_n(1;x)-1\|_{\psi_1}=\alpha_n,\eqno(3.7)
\end{equation*}
\begin{equation*}
    \|L_n(\rho;x)-\rho\|_{\psi_2}=\beta_n,\eqno(3.8)
\end{equation*}
\begin{equation*}
    \|L_n(\rho^2;x)-\rho^2\|_{\psi_3}=\gamma_n,\eqno(3.9)
\end{equation*}
where $\alpha_n,~\beta_n$ and $\gamma_n$ tend to zero as $n\to\infty$ and $\psi(x)=\max\big\{\psi_1(x),\psi_2(x),\psi_3(x)\big\}$. Then for all
$f\in C_\rho^k(\mathbb{R}^+_0)$ the inequality
\begin{equation*}
    \|L_n(f;x)-f(x)\|_{\psi\rho^2}\leq16\Omega_\rho\big(f,\sqrt{\alpha_n+2\beta_n+\gamma_n}\big)+\alpha_n\|f\|_\rho
\end{equation*}
holds for $n$ large enough.

\vspace{.25cm}
\parindent8mm Now we define following a sequence of positive linear operators $P_n^*$ with the help of $T_n^*$ defined by (1.7)

\begin{equation*}
P_n^*(f;x):=\frac{\rho^2(x)e^{-\frac{n}{b_n}xH(1)}}{A(1)}\sum\limits_{k=0}^{\infty}\frac{f\big(\frac knb_n\big)}{\rho^2\big(\frac knb_n\big)}p_k\bigg(\frac{n}{b_n}x\bigg).\eqno(3.10)
\end{equation*}

\vspace{.5cm}
\noindent\noindent\textbf{Theorem 3.7.} Let $P_n^*$ be the sequence of the positive linear operators defined by (3.10) and $\psi(x)=1+x^2$. If
$f\in C_\rho^k(\mathbb{R}^+_0)$, then
\begin{equation*}
    \|P_n^*(f;x)-f(x)\|_{\rho^4\psi}\leq16\Omega_\rho\big(f,\sqrt{\alpha_n+2\beta_n}\big)+\alpha_n\|f\|_\rho
\end{equation*}

\noindent\noindent\textbf{Proof.} By simple calculations we get
\begin{equation*}
    P_n^*(1;x)-1=\rho^2(x)\bigg[\frac{e^{-\frac{n}{b_n}xH(1)}}{A(1)}\sum\limits_{k=0}^{\infty}\frac{1}{\rho^2\big(\frac knb_n\big)}p_k\bigg(\frac{n}{b_n}x\bigg)-\frac1{\rho^2(x)}\bigg],\eqno(3.11)
\end{equation*}

\begin{equation*}
    P_n^*(\rho;x)-\rho(x)=\rho^2(x)\bigg[\frac{e^{-\frac{n}{b_n}xH(1)}}{A(1)}\sum\limits_{k=0}^{\infty}\frac{1}{\rho\big(\frac knb_n\big)}p_k\bigg(\frac{n}{b_n}x\bigg)-\frac1{\rho(x)}\bigg],\eqno(3.12)
\end{equation*}

\begin{equation*}
    P_n^*(\rho^2;x)-\rho^2(x)=0.\eqno(3.13)
\end{equation*}

\noindent From (3.3) and (3.5) we have
\begin{equation*}
    \lim\limits_{n\to\infty}\bigg\|\frac{e^{-\frac{n}{b_n}xH(1)}}{A(1)}\sum\limits_{k=0}^{\infty}\frac{1}{\rho^2\big(\frac knb_n\big)}p_k\bigg(\frac{n}{b_n}x\bigg)-\frac1{\rho^2(x)}\bigg\|_\psi=0,
\end{equation*}

\begin{equation*}
    \lim\limits_{n\to\infty}\bigg\|\frac{e^{-\frac{n}{b_n}xH(1)}}{A(1)}\sum\limits_{k=0}^{\infty}\frac{1}{\rho\big(\frac knb_n\big)}p_k\bigg(\frac{n}{b_n}x\bigg)-\frac1{\rho(x)}\bigg\|_\psi=0
\end{equation*}

\noindent using (3.4) and (3.11) we obtain

\begin{eqnarray*}
  \|P_n^*(1;x)-1\|_{\rho^2\psi}&=& \bigg\|\frac{e^{-\frac{n}{b_n}xH(1)}}{A(1)}\sum\limits_{k=0}^{\infty}\frac{1}{\rho^2\big(\frac knb_n\big)}p_k\bigg(\frac{n}{b_n}x\bigg)-\frac1{\rho^2(x)}\bigg\|_\psi  \\
  &\leq& \frac{b_n}{n}\frac{A(1)
+2A^\prime(1)+A(1)H^{\prime\prime}(1)}{A(1)}+\frac{b_n^2}{n^2}\frac{A^\prime(1)+A^{\prime\prime}(1)}{A(1)}\\
&=&\alpha_n.
\end{eqnarray*}

\noindent By means of (3.2) and (3.12), one gets

\begin{eqnarray*}
  \|P_n^*(\rho;x)-\rho\|_{\rho^2\psi}&=& \bigg\|\frac{e^{-\frac{n}{b_n}xH(1)}}{A(1)}\sum\limits_{k=0}^{\infty}\frac{1}{\rho\big(\frac knb_n\big)}p_k\bigg(\frac{n}{b_n}x\bigg)-\frac1{\rho(x)}\bigg\|_\psi  \\
  &\leq& \frac{b_n}{n}\frac{A^\prime(1)}{A(1)}\\
&=&\beta_n.
\end{eqnarray*}

\noindent finally from (3.13), we obtain

\begin{eqnarray*}
  \|P_n^*(\rho^2;x)-\rho^2\|_{\rho^2\psi}&=& 0\\
&=&\gamma_n.
\end{eqnarray*}

Thus the (3.7)-(3.9) assumptions of Theorem 3.6 are satisfied for the operators (3.10). From Theorem 3.6, we have

\begin{equation*}
    \|P_n^*(f;x)-f(x)\|_{\rho^4\psi}\leq16\Omega_\rho\big(f,\sqrt{\alpha_n+2\beta_n+\gamma_n}\big)+\alpha_n\|f\|_\rho
\end{equation*}

\noindent for each $f\in C_\rho^k(\mathbb{R}_0^+)$. This completes the proof.

\vspace{.5cm}
\noindent\textbf{4. Numerical Examples}\\

\vspace{.25cm}
\noindent\textbf{Example 4.1.} The sequence $\{(1+x)^k\}_{k=1}^\infty$ which is the Sheffer sequence for $A(t)=e^t$ and $H(t)=t$ has the generating function of the following type
\begin{equation*}
  e^{(1+x)t}=\sum\limits_{k=0}^{\infty}\frac{(1+x)^k}{k!}t^k.
\end{equation*}
Let us select $p_k(x)=\frac{(1+x)^k}{k!}$. Making use of above knowledge $p_k(x)\ge0$ for $x\in[0,\infty),~A(1)\ne0$ and $H^\prime(1)=1$ are provided. Considering these polynomials in (1.7), we obtain operators as follows
\begin{equation*}
  T_n^*(f;x)=e^{-(\frac n{b_n}x+1)}\sum\limits_{k=0}^{\infty}\frac{(\frac n{b_n}x+1)^k}{k!}f\bigg(\frac kn b_n\bigg).
\end{equation*}
The error bound for the function $f(x)=-4xe^{-3x}$ under the condition condition $A(t)=e^t$ and $H(t)=t$ is computed in the following Table 1:
\begin{table}[ h! ]
\centering
\begin{tabular}{ p{.5cm} p{11.5cm}}
\hline
\hline
n & Error estimate by $T_n^*$ operators including $\{(1+x)^k\}_{k=1}^\infty$ sequence\\
\hline
10 & \hspace{4cm} 0.9481710727 \\
$10^3$ & \hspace{4cm} 0.8474426939 \\
$10^5$ & \hspace{4cm} 0.3806348279 \\
$10^7$ & \hspace{4cm} 0.1348930985 \\
$10^9$ & \hspace{4cm} 0.0442354247 \\
$10^{11}$ & \hspace{4cm} 0.0141505650 \\
$10^{13}$ & \hspace{4cm} 0.0044911482 \\
$10^{15}$ & \hspace{4cm} 0.0014218648 \\
$10^{17}$ & \hspace{4cm} 4.4979717260e-004 \\
$10^{19}$ & \hspace{4cm} 1.4225476356e-004 \\
\hline
\end{tabular}
\caption{The error bound of function $f(x)=-4xe^{-3x}$ by using
modulus of continuity}
\label{table:1}
\end{table}


\vspace{1cm}
\noindent\textbf{Example 4.2.} The sequence $\{x^{k-1}\}_{k=1}^\infty$ which is the Sheffer sequence for $A(t)=t$ and $H(t)=t$ has the generating function of the following type
\begin{equation*}
  te^{xt}=\sum\limits_{k=1}^{\infty}\frac{x^{k-1}}{(k-1)!}t^k.
\end{equation*}
Let us select $p_k(x)=\frac{x^{k-1}}{(k-1)!}$. Making use of above knowledge $p_k(x)\ge0$ for $x\in[0,\infty),~A(1)\ne0$ and $H^\prime(1)=1$ are provided. Considering these polynomials in (1.7), we obtain operators as follows
\begin{equation*}
  T_n^*(f;x)=e^{-\frac n{b_n}x}\sum\limits_{k=1}^{\infty}\frac{(\frac n{b_n}x)^{k-1}}{(k-1)!}f\bigg(\frac kn b_n\bigg).
\end{equation*}
The error bound for the function $f(x)=-4xe^{-3x}$ under the condition condition $A(t)=t$ and $H(t)=t$ is computed in the following Table 2:
\begin{table}[ h! ]
\centering
\begin{tabular}{ p{.5cm} p{11.5cm}}
\hline
\hline
n & Error estimate by $T_n^*$ operators including $\{x^{k-1}\}_{k=1}^\infty$ sequence\\
\hline
10 & \hspace{4cm} 0.8938844531 \\
$10^3$ & \hspace{4cm} 0.8409966996 \\
$10^5$ & \hspace{4cm} 0.3803350939 \\
$10^7$ & \hspace{4cm} 0.1348824385 \\
$10^9$ & \hspace{4cm} 0.0442350750 \\
$10^{11}$ & \hspace{4cm} 0.0141505538 \\
$10^{13}$ & \hspace{4cm} 0.0044911478 \\
$10^{15}$ & \hspace{4cm} 0.0014218647 \\
$10^{17}$ & \hspace{4cm} 4.4979717224e-004 \\
$10^{19}$ & \hspace{4cm} 1.4225476355e-004 \\
\hline
\end{tabular}
\caption{The error bound of function $f(x)=-4xe^{-3x}$ by using
modulus of continuity}
\label{table:2}
\end{table}

\vspace{.5cm}
\noindent\textbf{Conclusion}\newline

\noindent We introduced the Chlodowsky variant of generalized Sz\`{a}sz operators by means of Sheffer polynomials and established different approximation results. We have also given an algorithm to plot the graphs of the positive linear operators and with the help of these graphical examples, we claimed that our operators are better than the old operators to approximate a given function. Some numerical examples are also provided and we found the error bound of a given function by using modulus of smoothness.

\end{document}